\documentstyle[12pt]{article}

\titlepage

\title{Holonomy and four-dimensional manifolds}
\author {Richard Atkins}
\date{}
\newtheorem{fact}{Fact}

\newtheorem{lemma}[fact]{Lemma}

\newtheorem{theorem}[fact]{Theorem}

\begin{document}
\maketitle
Department of Mathematics, Trinity Western University, Glover Road, Langley, 
British Columbia V2Y 1Y1, Canada \\ 
\section{Introduction}
This paper investigates the relationship between two fundamental types of objects 
associated with a connection on a manifold: the existence of parallel semi-Riemannian 
metrics and the associated holonomy group. Typically in Riemannian geometry, 
a metric is specified which determines a Levi-Civita connection. Here we consider 
the connection as more fundamental and allow for the possibility of several parallel
metrics. Holonomy is an old geometric concept which is enjoying revived interest in 
certain branches of mathematical physics, in particular loop quantum gravity and 
Calabi-Yau manifolds in string theory. It measures, in group theoretic terms the connection's deviation 
from flatness and takes the topology of the manifold into account. 

It is well known that for a Riemannian manifold, the reducibility of the holonomy
group of the Levi-Civita  connection implies the existence of multiple independent
parallel Riemannian metrics on the manifold \cite{aa}. In this paper we 
look at the converse: when does the existence of multiple independent semi-Riemannian
metrics on a manifold, parallel with respect to a linear connection, imply the 
reducibility of the holonomy group of the connection? We do not necessarily
assume that the connection is symmetric in our solutions to this problem.
If $g$ is a semi-Riemannian metric on $M$ which is parallel with respect to the 
connection $\nabla$, then the holonomy group $\Psi(x)$, at $x\in M$, preserves 
$g$: $\psi^{*}(g)=g$, for all $\psi\in \Psi(x)$. Thus the existence of parallel metrics
places algebraic restrictions on $\Psi(x)$; these restrictions will be  the subject 
of our investigations. For manifolds of dimension $d\neq 4$ the problem has a 
purely algebraic solution. For four-dimensional manifolds the relationship
of the parallel metrics 
of the connection to the reducibility of the holonomy group is not entirely
algebraic but depends also on the fundamental group of the manifold: the presence of 
three parallel metrics always implies reducibility but reducibiliy in the case of two 
metrics of signature $(2,2)$ is guaranteed only for simply connected manifolds. The central 
theorem in this paper is the construction of a topologically non-trivial four-dimensional 
manifold with a connection that admits two independent metrics of signature $(2,2)$ and yet 
has irreducible holonomy. $d=4$ is the critical dimension with respect to reducibility of 
the holonomy group.

It is interesting to note that $d=4$ appears as the critical dimension in other conttexts as well. In quantum
field theory, for instance, infinite divergences appear in the calculation of 
scattering amplitudes as the dimension of spacetime approaches four. Also, it has been shown that 
$R^{4}$ has the remarkable property of admitting exotic differentiable structures. 
In superstring theory spacetime is ten or eleven dimensional but only four dimensions are observed in nature.
Therefore some unique characteristic of four-dimensional manifolds must be involved in explaining
this mismatch of dimensions.
 
The problem of non-uniqueness of parallel metrics, largely in the case of
Lorentzian connections,  has been investigated by several authors 
 (cf. \cite{bb}-\cite{ff}, \cite{gg}-\cite{ii}).  

\section{Reducibility of the Holonomy Group}

\label{Reducibility of the Holonomy Group}

Let $V$ be a vector space over a field ${F}$, which will be either ${R}$
or ${C}$.  Let $G$ be a group which acts on $V$ on the left. A subspace
$W$ of $V$ is said to be $G$-{\it invariant} if $g\cdot \xi \in  W$, for all
$g\in G$ and $\xi \in W$. If there exists a proper, non-trivial $G$-invariant
subspace of $V$ we say that $G$ {\it acts reducibly}  on $V$, or more simply,
that $G$ is {\it reducible}. In our applications,
$G$ shall be a subgroup of $Aut(V)$, the group of linear automorphisms of $V$.
The holonomy group $\Psi(x)$ of a linear connection $\nabla$ on a connected
manifold $M$, at $x\in M$, is a subgroup of $Aut (T_{x}M)$. For any two points
$x,y\in M$, the holonomy groups $\Psi(x)$ and $\Psi(y)$ are isomorphic,
since $M$ is connected. If for some point (and hence all points) $x\in M$, $\Psi(x)$
is reducible (respectively, irreducible) then we say that the connection $\nabla$
has  {\it reducible} (respectively, {\it irreducible}) {\it holonomy}. 
The results of this section are essentially algebraic and may be derived with the aid
of standard normal form theorems in linear algrbra. Since the proofs are somewhat 
lengthy and
technical they shall be omitted.

We begin with a theorem that provides sufficient conditions, 
with regard to the existence of parallel semi-Riemannian metrics, to ensure the 
reducibility of the holonomy group of the connection. 

\begin{theorem}
\hspace{3in}

(i) Let $g_{1}$ and $g_{2}$ be two independent semi-Riemannian metrics on a connected
manifold $M$, parallel with respect to a linear connection $\nabla$ on $M$. Suppose
that one of $g_{1}$, $g_{2}$ has signature $(p,q)$ with $p\neq q$. Then
$\nabla$ has reducible holonomy. 

(ii) Let $g_{1}$ and $g_{2}$ be two independent semi-Riemannian metrics on a
connected, two-dimensional manifold $M$, parallel with respect to a linear connection 
$\nabla$ on $M$. Then $\nabla$ has reducible holonomy.

(iii) Suppose $n\not\equiv 0$ $mod$ $4$.
Let $g_{1}$, $g_{2}$ and $g_{3}$ be three 
independent semi-Riemannian metrics on a connected, $n$-dimensional
manifold $M$, parallel with respect 
to a linear connection $\nabla$ on $M$. Then $\nabla$ has reducible holonomy.

(iv) Let $g_{1}$, $g_{2}$, $g_{3}$ and $g_{4}$ be four independent, semi-Riemannian 
metrics on a connected manifold $M$, parallel with respect to a linear connection
$\nabla$ on $M$. Then $\nabla$ has reducible holonomy.
\label{theorem:reduce2}
\end{theorem}

It is possible to construct examples of connections on manifolds that show 
that the numbers of parallel semi-Riemannian metrics in the above theorem are sharp. 
Specifically we have the following.

\begin{theorem}
\hspace{3in}

(i) Let $n=2m$ and $m\geq 3$. There exist two independent semi-Riemannian
metrics on ${R}^{n}$ of signature $(m,m)$, parallel with respect to a linear
connection having an irreducible holonomy group. 

(ii) Let $n=2m$, $m=2r$ and $r\geq 2$.
There exist three independent semi-Riemannian
metrics on ${R}^{n}$ of signature $(m,m)$, parallel with respect to a linear
connection having an irreducible holonomy group.
\label{theorem:irreduce2}
\end{theorem}

\section{Four-Dimensional Manifolds}
\label{Four-Dimensional Manifolds}

\subsection{Simply connected four-dimensional manifolds}

The following theorem is also proved algebraically and uses the fact that 
the holonomy group of a connection on a simply connected manifold consists
of only one component (cf. \cite{aa}, p 73).

\begin{theorem}
\hspace{3in}

(i) Let $M$ be a simply connected, four-dimensional manifold.
Let $g_{1}$ and $g_{2}$ be two independent semi-Riemannian metrics
on $M$, parallel with respect to a linear connection $\nabla$ on $M$. 
Then $\nabla$ has  reducible holonomy.

(ii) Let $M$ be a connected, four-dimensional manifold. 
Let $g_{1},g_{2}$ and $g_{3}$ be three independent semi-Riemannian
metrics on $M$, parallel with respect to a linear connection $\nabla$ on $M$. Then
$\nabla$ has reducible holonomy.
\label{theorem:4dim3}
\end{theorem}

\subsection{Irreducibility of the holonomy group}

In this section we shall construct an example of a (non-simply connected)
four-dimensional manifold endowed with two independent semi-Riemannian 
metrics of signature $(2,2)$ which are parallel with respect to a 
linear connection  having an irreducible holonomy group.

Let ${\cal {F}}$ be the subgroup of $GL(4;{R})$ whose elements $F$
are of the form
\[ F = \left( \begin{array}{rr}
A & B \\
-B & A \end{array} \right).\]
The map $\phi:{\cal {F}} \rightarrow GL(2;{C})$ defined by
\[\phi(F) := A+iB \]
is a group monomorphism.  Let ${\cal {G}}$ be any subgroup of ${\cal {F}}$.
We denote by ${\cal {G}}_{{C}}$ the subgroup $\phi({\cal {G}})$ of
$GL(2;{C})$.  $\phi:{\cal {G}}\rightarrow {\cal {G}}_{{C}}$
is a Lie group isomorphism. Let  ${C}_{{R}}^{2}$ denote the vector space  of 
$2$-tuples of complex
numbers over the real field. The following lemma is easily established.

\begin{lemma}
${\cal {G}}$ acts reducibly on ${R}^{4}$ if and only if
${\cal {G}}_{{C}}$ acts reducibly on ${C}_{{R}}^{2}$.

\label{lemma:complex2}

\end{lemma}

Let $K$ and $L$ denote the matrices in $GL(4;{R})$ defined as follows
\[ K := \left( \begin{array}{cc} I_{2\times 2} & 0 \\ 0 & -I_{2\times 2} \end{array}
\right) \hspace*{.4in} and  \hspace*{.4in}
 L := \left( \begin{array}{cc} 0 &  I_{2\times 2} \\ I_{2\times 2} & 0 \end{array}
\right). \]
$K$ and $L$ may be regarded as symmetric bilinear forms on ${R}^{n}$ with 
signature $(2,2)$.
Let ${\cal {H}}$ be the subgroup of $GL(4;{R})$ whose elements $H$
satisfy
\[ \left\{ \begin{array}{lll}
^{t}HKH & = & K, \hspace*{.1in} and \\
^{t}HLH & = & L. \end{array}  \right. \]
Let $G \in \cal{F}$ be defined by
\[ G:= \left( \begin{array}{cc|cc}
 \sqrt{2} & 0 & 0 & 1 \\
 0 & -\sqrt{2} & 1 & 0 \\ \hline 
 0 & -1 & \sqrt{2} & 0   \\ 
 -1 & 0 & 0 & -\sqrt{2}   \end{array} \right), \]
Then $G \in \cal{H}$ and $G^2=I$.
Let ${\cal {G}}$ denote the subgroup  of $\cal{F}$ consisting of
the two elements $I$ and $G$: ${\cal {G}}:= \{I,G\}$. The characteristic 
polynomial $p$ of $G$ is given by $p(x)=(x^{2}-1)^{2}$. The eigenspace
$E$ of $G$ corresponding to the eigenvalue $\lambda=1$ is a two-dimensional
subspace of ${R}^{4}$. Set $M':={R}^{4}-E$. ${\cal {G}}$ acts 
properly discontinously (cf. \cite{aa} Vol I, pp 43-44) 
on the left of $M'$, by  matrix multiplication.
Let $g_{1}$ and $g_{2}$ denote the 
semi-Riemannian metrics on $M'$ represented, respectively,  by $K$ and $L$ 
in the moving frame $(dx^{1},...,dx^{4})$ on $M'$. 
We see that $L_{G}^{*}(g_{i}) = g_{i}$, for $i=1,2$, where $L_{G}$ denotes 
left multiplication by $G$.

Define the  1-forms $\alpha$ and $\beta$ on $M'$ by
\[ \alpha := x^{1}dx^{4}-x^{4}dx^{1} \hspace*{.3in} and  \hspace*{.3in}
\beta := x^{2}dx^{3}-x^{3}dx^{2}\]
and the matrix of 1-forms, $\theta$, by
\[ \theta:= \left( \begin{array}{cc|cc}
0 & \alpha & 0 & \beta \\
-\alpha & 0 & -\beta & 0 \\ \hline
0 & -\beta & 0 & \alpha \\
\beta & 0 & -\alpha & 0  \end{array} \right). \]
The following two lemmas may be demonstrated by direct calculation.

\begin{lemma}
$G\theta =-\theta G .$

\label{lemma:2forms1}

\end{lemma}

\begin{lemma}
$L_{G}^{*}(\theta) = -\theta.$

\label{lemma:2forms2}

\end{lemma}

Define the linear connection $\nabla$ on $M'$ by
\[ \nabla_{X_{i}}dx^{j} := 
-\sum_{l=1}^{4}\theta^{j}_{l} (X_{i})dx^{l} ,\]
for $1 \leq i,j\leq 4$, where $X_{i}=\partial/\partial x^{i}$.

\begin{lemma}
$g_{1}$ and $g_{2}$ are parallel with respect to $\nabla$.

\label{lemma:2forms3}

\end{lemma}
{\bf Proof:}
A metric $g=\sum_{i,j=1}^{4}C_{ij}dx^{i}\otimes dx^{j}$ on $M'$, 
with $C \in GL(4;{R})$,  is parallel with respect to $\nabla$
if and only if $^{t}\theta C + C \theta =0$. The lemma follows from
the fact that $^{t}\theta K + K \theta =0$ and 
$^{t}\theta L + L \theta =0$. \\
{\bf Q.E.D.} \\

For a vector field $X$ on $M'$, we denote $L_{G}(X):= L_{G*}\circ X\circ L_{G}$.

\begin{lemma}
Let $X,Y$ be vector fields on $M'$. Then
\[ L_{G}(\nabla_{X}Y) = \nabla_{L_{G}(X)}L_{G}(Y).\]

\label{lemma:2forms4}

\end{lemma}
{\bf Proof:}
It suffices to consider $X=X_{i}=\partial/\partial x^{i} $ and 
$Y=X_{j}=\partial/\partial x^{j} $.
\begin{eqnarray*}
L_{G}(\nabla_{X_{i}}X_{j})(x)
& = \hspace*{.1in} & 
\sum_{s=1}^{4}\sum_{t=1}^{4}G^{t}_{s}\theta^{s}_{j}(X_{i}|_{Gx})X_{t}(x) \\
& = \hspace*{.1in} &  \hspace*{-.15in}
-\sum_{s=1}^{4}\sum_{t=1}^{4}G^{s}_{j}\theta^{t}_{s}(X_{i}|_{Gx})X_{t}(x), 
\hspace*{.62in} by \hspace*{.1in} Lemma \hspace*{.1in} \ref{lemma:2forms1}\\
& = \hspace*{.1in} & 
\sum_{s=1}^{4}\sum_{t=1}^{4}G^{s}_{j}L_{G}^{*}(\theta^{t}_{s})(X_{i}|_{Gx})X_{t}(x),
\hspace*{.3in} by \hspace*{.1in} Lemma \hspace*{.1in} \ref{lemma:2forms2} \\
& = \hspace*{.1in} & \nabla_{L_{G}(X_{i})}L_{G}(X_{j})(x). 
\end{eqnarray*}
{\bf Q.E.D.}\\

In summary, the left group action of ${\cal {G}}$ on $M'$ preserves
$g_{1}$, $g_{2}$ and $\nabla$. It follows that $g_{1}$, $g_{2}$ and $\nabla$
project to the quotient manifold
\[ M:=M'/{\cal {G}}.\]
We denote the projections of $g_{1}$, $g_{2}$ and $\nabla$ to $M$ by the
same symbols.

Now let $\gamma :[0,1]\rightarrow M'$ be the smooth curve defined by
\[ \gamma (t) := \left( \begin{array}{c}
1-2t \\ t(t-1) \\ t(t-1) \\ (\sqrt{2}+1)(2t-1) \end{array} \right).\]
Set $x_{0}:=\gamma (0)$ and $x_{1}:= \gamma (1)$. Since $x_{1}=Gx_{0}$,
$\gamma$ is a closed loop in $M$. Let $\tau':T_{x_{0}}M' \rightarrow T_{x_{1}}M'$ denote 
parallel translation along $\gamma$ in $M'$ and let 
$\tau:T_{x_{0}}M \rightarrow T_{x_{0}}M$ denote parallel 
translation along $\gamma$ in $M$. Let $a^{j}X_{j}|_{x_{0}}\in T_{x_{0}}M'$, 
where we have used the summation convention. 
Then $\tau'(a^{j}X_{j}|_{x_{0}})= f^{j}(1)X_{j}|_{x_{1}}$, where 
$f:[0,1]\rightarrow {R}^{4}$ is the unique curve satisfying
\[ \left\{ \begin{array}{llll}
\dot{f}+\theta (\dot{\gamma})f &  = & 0, & \hspace*{.1in} and \\
f^{j}(0) & = & a^{j}, & \hspace*{.1in} 1\leq j\leq 4. 
\end{array} \right. \] 
Now $\theta (\dot{\gamma}) \equiv 0$ and so $f = f(t)$ is constant.
Therefore $\tau'(a^{j}X_{j}|_{x_{0}}) = a^{j}X_{j}|_{x_{1}}.$
It follows that 
\[ \tau(a^{j}X_{j}|_{x_{0}}) = G_{k}^{j}a^{k}X_{j}|_{x_{0}}.\]
We identify $T_{x_{0}}M$ with ${R}^{4}$ by means of the basis 
$(X_{1},...,X_{4})$ and hence the holonomy group $\Psi (x_{0})$ of $\nabla$
at $x_{0}\in M$ is identified with elements of $GL(4;{R})$. 
We have shown above that $G\in \Psi (x_{0})$.

The curvature form $\Omega$ of $\nabla$ is 
\[ \Omega  =   d\theta +\theta\wedge\theta = d\theta =  
\left( \begin{array}{cc|cc}
0 & d\alpha & 0 & d\beta \\
-d\alpha & 0 & -d\beta & 0 \\ \hline
0 & -d\beta & 0 & d\alpha \\
d\beta & 0 & -d\alpha & 0  \end{array} \right). \]
Now $d\alpha = 2dx^{1}\wedge dx^{4}$ and $d\beta = 2dx^{2}\wedge dx^{3}$.  
The Lie algebra $\psi^{inf}(x_{0})$ of the infinitesimal holonomy group 
$\Psi^{inf}(x_{0})$ of $\nabla$ at $x_{0}\in M$ consists of elements 
$S\in \cal{F}$ of the form
\[ S:= \left( \begin{array}{cc|cc}
0 & a & 0 & b \\
-a & 0 & -b & 0 \\ \hline
0 & -b & 0 & a \\
b & 0 & -a & 0  \end{array} \right), \]
where $a,b\in {R}$.
The Lie algebra $\psi^{inf}_{{C}}(x_{0})$ of the complexified infinitesimal 
holonomy group $\Psi^{inf}_{{C}}(x_{0})$ consists of elements 
$s\in GL(2;C)$ of the form 
\[ s = \left( \begin{array}{cc}
0 & a \\ -a & 0 \end{array}  \right), \]
where $a\in {C}$. That is, $\psi^{inf}_{{C}}(x_{0})= so(2;{C})$ 
and so $\Psi^{inf}_{{C}}(x_{0}) = SO(2;{C})$.

Let ${\cal {A}}$ denote the $R-$subalgebra of $GL(2;{C})$ generated by
$g$ and $SO(2;{C})$ over $R$. It is not difficult to demonstrate the
following lemma.

\begin{lemma}
${\cal {A}}$ = $gl(2;{C})$.

\label{lemma:2forms5}

\end{lemma}

Let ${\cal {H}}$ be the subgroup of $GL(4;{R})$ such that
${\cal {H}}_{C}$ is the subgroup of $GL(2;{C})$  generated by $g$
and $SO(2;{C})$. By Lemma \ref{lemma:2forms5}, ${\cal {H}}_{C}$ acts 
irreducibly on ${C}^{2}_{R}$. By Lemma \ref{lemma:complex2},
${\cal {H}}$ acts irreducibly on ${R}^{4}$. Since ${\cal {H}}\subseteq\Psi (x_{0})$, 
$\Psi (x_{0})$ acts irreducibly on $T_{x_{0}}M$. We arrive at the
following theorem.

\begin{theorem}
$g_{1}$ and $g_{2}$ are  two independent semi-Riemannian metrics of signature
$(2,2)$ on a non-simply connected, four-dimensional manifold $M$, parallel with
respect to a linear connection $\nabla$ on $M$. Moreover, $\nabla$ has
irreducible holonomy.		

\label{theorem:2forms1}

\end{theorem}

\end{document}